\documentclass[12pt]{article}
\usepackage{amssymb,amsmath,epsfig,amscd,eucal,psfrag}

\pdfoutput=1

\newtheorem{thm}{Theorem}
\newtheorem{defn}[thm]{Definition}
\newtheorem{lem}[thm]{Lemma}

\newtheorem{cor}[thm]{Corollary}

\newtheorem{qu}[thm]{Question}
\newtheorem{letterthm}{Theorem}
\newtheorem{lettercor}[letterthm]{Corollary}

\newenvironment{pf}{\par\medskip\noindent{\em Proof. }}{\hfill $\square$\par\medskip}

\newcommand{\C}{\mathbb{C}}

\newcommand{\F}{\mathbb{F}}

\newcommand{\Q}{\mathbb{Q}}

\newcommand{\Z}{\mathbb{Z}}

\newcommand{\into}{\hookrightarrow}

\input xy
\xyoption{all}

\title{On surface subgroups of doubles of free groups}
\author{Cameron Gordon and Henry Wilton}

\begin{document}

\maketitle

\begin{abstract}
We give several sufficient conditions for a double of a free group along a cyclic subgroup to contain a surface subgroup.
\end{abstract}

\section{Introduction}

\begin{qu}[Gromov \cite{bridson_problems_????}]\label{q: Gromov's Question}
Does every one-ended (word-)hyperbolic group contain a surface subgroup?
\end{qu}

By a `surface subgroup' we mean a subgroup isomorphic to the fundamental group of a closed surface with non-positive Euler characteristic.  Of course, in a hyperbolic group such a surface must actually be of negative Euler characteristic.  However, the techniques of this paper apply equally well to some non-hyperbolic groups.

Very little is known about Gromov's question even for some very concrete classes of hyperbolic groups.  For example, if $F_n$ is the free group of rank $n$ then it follows from Bestvina and Feighn's combination theorem \cite{bestvina_combination_1992} that the double
\[
D_n(w)=F_n*_{\langle w\rangle}F_n
\]
is word-hyperbolic if and only if $w\in F_n$ is not a proper power.  Even in this class of examples, the answer to Gromov's question remains unknown.  It has been suggested that Gromov's question may have a negative answer for doubles.

The following recent theorem of Calegari represents the first real progress.  Throughout, we shall use $\beta_i$ to denote the $i$th Betti number of a group or a topological space.

\begin{thm}[Calegari \cite{calegari_surface_2008}]\label{t: Calegari's Theorem}
If a hyperbolic group $G$ is the fundamental group of a graph of free groups with cyclic edge groups and $\beta_2(G)>0$ then $G$ contains a surface subgroup.
\end{thm}

Calegari's theorem reduces Gromov's question for the double $D_n(w)$ to a condition on the second virtual Betti number.   For a group $G$, the \emph{$i$th virtual Betti number} is the supremum of $\beta_i(G')$ where $G'$ ranges over all finite-index subgroups of $G$.

\begin{cor}
Let $w\in F_n$.  The virtual second Betti number of $D_n(w)$ is positive if and only if $D_n(w)$ contains a surface subgroup.
\end{cor}
\begin{pf}
If $w=u^r$ is a proper power then $D_n(w)$ both contains a subgroup isomorphic to $\Z^2$ and has positive second virtual Betti number.  We sketch the proofs of these facts.  Let $u_1$ and $u_2$ be the copies of $u$ in the two vertex groups of $D_n(w)$.  Let $v=u_1^{-1}u_2$.  Clearly $v$ and $w$ commute.  By Marshall Hall's Theorem \cite{hall_subgroups_1949}, $F_n$ has a finite-index subgroup $F'$ such that $w$ is an element of a basis for $F'$.  Let $\xi:D_n(w)\to F_n$ be the natural map given by identifying the two vertex groups and let $D'=\xi^{-1}(F')$.  The inclusion map $\langle w\rangle\into D'$ has a left inverse $\lambda$ (that is, $\lambda$ is a \emph{retraction}).  Note that $\lambda(v)=1$.  Also, $v$ maps to a simple loop in the underlying graph of the induced graph-of-groups decomposition of $D'$, and so there is a retraction $\mu:D'\to\langle v\rangle$ such that $\mu(w)=1$.  Therefore, the map $(\lambda,\mu)$ is a retraction $D'\to\langle v,w\rangle\cong\Z^2$.   It follows that the induced map $H_*(\Z^2)\to H_*(D')$ is injective.  In particular, $\beta_2(D')$ is positive.

Suppose therefore that $w$ is not a proper power so that, as remarked above, $D_n(w)$ is hyperbolic \cite{bestvina_combination_1992}.   If $D_n(w)$ has a finite-index subgroup $D'$ with $\beta_2(D')>0$ then it is an immediate consequence of Theorem \ref{t: Calegari's Theorem} that $D'$, and hence $D$, contains a surface subgroup.  Conversely, suppose that $D_n(w)$ contains a surface subgroup $S$.  Because $w$ has no proper roots, $D_n(w)$ is a limit group.  By Theorem B of \cite{wilton_halls_2008}, $S$ is a \emph{virtual retract} of $D_n(w)$---that is, $D_n(w)$ has a finite-index subgroup $D'$ that contains $S$ and such that the inclusion map $S\hookrightarrow D'$ has a left inverse.  As above, $H_*(S)$ embeds into $H_*(D')$ and so $\beta_2(D')>0$.
\end{pf}

The following question is therefore equivalent to Question \ref{q: Gromov's Question} for hyperbolic doubles.

\begin{qu}\label{q: b_2}
Does every one-ended double $D_n(w)$ have a finite-index subgroup $D'$ with $\beta_2(D')>0$?
\end{qu}

One can show that $D_n(w)$ is one-ended if and only if $w$ is not contained in a proper free factor of $F_n$.  We leave this as an exercise to the reader.

Nothing is known about Gromov's question for groups of the form $D_n(w)$ in the absence of Calegari's hypothesis of positive second Betti number, which is equivalent, in the case of doubles, to the condition that $w\in [F_n,F_n]$.  In this paper we present various approaches to Gromov's question, and provide several infinite families of new examples of doubles with surface subgroups.

As well as determining a double, a word in a free group also determines a one-relator group.  For $w\in F_n$, let
\[
G_n(w)=F_n/\langle\langle w\rangle\rangle
\]
the associated one-relator group.  There is a homomorphism $\zeta: D_n(w)\to G_n(w)$ obtained by concatenating the natural maps $\xi:D_n(w)\to F_n$ and $\eta:F_n\to G_n(w)$.  Our first theorem relates the homology of subgroups of $G_n(w)$ to the homology of subgroups of $D_n(w)$.  Throughout this paper, homology is assumed to have integral coefficients, unless otherwise indicated.

\begin{letterthm}\label{t: Theorem A}
Let $w\in F_n$.  Let $G'$ be a subgroup of $G_n(w)$ and let $D'=\zeta^{-1}(G')$.  There is a surjection
\[
\psi:H_2(D')\to H_2(G').
\]
Furthermore, if $w$ is not a proper power then $\psi$ is an isomorphism.
\end{letterthm}

A little care is needed here---$\psi$ is not induced by $\zeta$.  Indeed, $\zeta$ factors through a free group and so is trivial at the level of second homology.

Theorem \ref{t: Theorem A} transfers the subject of Question \ref{q: b_2} from doubles to one-relator groups.  One-relator groups are often much more badly behaved than their corresponding doubles, so this may not be a great improvement.  (Indeed, every double $D_n(w)$ is a one-relator group!)  However, in the special case $n=2$, a one-relator group has zero Euler characteristic, which leads to a very simple relationship between the first and second virtual Betti numbers.

\begin{lettercor}\label{c: Corollary B}
Let $w\in F_n$.  If $G_n(w)$ has an index-$k$ subgroup $G'$ with
\[
\beta_1(G')>1+k(n-2)
\]
then $D_n(w)$ contains a surface subgroup.  In particular, in the case $n=2$, if $G_2(w)$ has a finite-index subgroup $G'$ with $\beta_1(G')> 1$ then $D_2(w)$ contains a surface subgroup.
\end{lettercor}

The study of one-relator groups has a long history, and one can use Corollary \ref{c: Corollary B} to provide large explicit families of new examples for which Gromov's question has a positive answer.  Any epimorphism $\phi:G_n(w)\to\Z$ provides a family of finite-index subgroups: $G_k=\phi^{-1}(k\Z)$.  The first homology of these covers is governed by the \emph{Alexander polynomial} $\Delta_\phi(t)$.

\begin{letterthm}[Theorem \ref{t: Roots of unity}]\label{t: Theorem D}
Let $w\in F_n$.  If there is an epimorphism $\phi:G_n(w)\to\Z$ such that $\Delta_\phi(t)$ has a root $\omega$ with $\omega^k=1$ for some $k$ then $\beta_1(G_k)>1+k(n-2)$.  Therefore $D_n(w)$ contains a surface subgroup.
\end{letterthm}

The Alexander polynomial is easy to compute, and every integral polynomial arises as an Alexander polynomial, so this is a large family of examples.  Perhaps more surprisingly, one can also gain information from the mod $p$ Alexander polynomial, using work of Howie \cite{howie_free_1998}.  The $n=2$ case is particularly simple.  In light of Theorem \ref{t: Calegari's Theorem}, we may assume that $w\in F_2\smallsetminus [F_2,F_2]$.  In this case there is a unique choice of $\phi$, and we write $\Delta_w(t)=\Delta_\phi(t)$.  Recall that a group is \emph{large} if a finite-index subgroup surjects a non-abelian free group.  A large group has virtually infinite first Betti number.

\begin{letterthm}[Theorem \ref{t: Alexander polynomial mod p}]\label{t: Theorem E}
Let $w\in F_2\smallsetminus [F_2,F_2]$.  If $\Delta_w(t)\equiv 0$ mod $p$ for some prime $p$ then $G_2(w)$ is large.  Therefore, $D_2(w)$ contains a surface subgroup.
\end{letterthm}

Baumslag and Pride proved that $G_n(w)$ is large whenever $n>2$ \cite{baumslag_groups_1978}.  We use results of Wise to find a family of large one-relator groups that includes two-generator examples.

\begin{letterthm}\label{t: Theorem C}
If $w\in F_n$ is a positive $C'(1/6)$ word then $G_n(w)$ is large.  Therefore, if $n=2$ then $D_2(w)$ contains a surface subgroup.
\end{letterthm}

Wise showed that one-relator groups are, in a suitable sense, generically $C'(1/6)$ (\cite{wise_residual_2001}, Theorem 6.1).  More generally, one can ask how common it is for a two-generator, one-relator group to be large or to have virtual first Betti number greater than one.   Button has used the Alexander polynomial to study large one-relator groups \cite{button_largeness_2008, button_proving_2008}.  In particular, he has shown that the vast majority of two-generator, one-relator presentations in `Magnus form' with cyclically reduced relation of length at most 12 are large (\cite{button_proving_2008}, Theorem 3.3).

However, there are examples of two-generator, one-relator groups with virtual first Betti number equal to one.  The most well-known such examples are Baumslag--Solitar groups, of the form $\langle a,b\mid b^{-1}a^pba^q\rangle$ for certain integers $p,q$ (see Theorem \ref{t: Edjvet-Pride} for the exact conditions on $p$ and $q$). We go on to find surface subgroups in the corresponding doubles using an entirely different method.

\begin{letterthm}\label{t: Theorem F}
Suppose $w=b^{-1}a^pba^q\in\langle a,b\rangle=F_2$.  Then $D_2(w)$ has a finite-index subgroup that is the fundamental group of a compact 3-manifold, and contains a surface subgroup.
\end{letterthm}

This paper is organized as follows.  In Section \ref{s: Mayer-Vietoris} we prove Theorem \ref{t: Theorem A} and deduce Corollary \ref{c: Corollary B}. We go on to characterize the circumstances under which Corollary \ref{c: Corollary B} can be expected to answer Question \ref{q: b_2}.  In Sections \ref{s: Alexander polynomial} and \ref{s: C'(1/6)} we apply Corollary \ref{c: Corollary B} to prove Theorems \ref{t: Theorem D}, \ref{t: Theorem E} and \ref{t: Theorem C}.  In Section \ref{s: Baumslag-Solitar groups} we give a new proof of Edjvet and Pride's theorem that Baumslag--Solitar groups have virtual first Betti number equal to one.  Finally, in Section \ref{s: Virtually geometric} we introduce a different approach to Question \ref{q: Gromov's Question} for doubles and prove Theorem \ref{t: Theorem F}.  We finish by asking whether this second approach applies to all doubles.

\section{A Mayer--Vietoris argument}\label{s: Mayer-Vietoris}

In this section we prove Theorem \ref{t: Theorem A} and deduce Corollary \ref{c: Corollary B}.  We will write $D=D_n(w)$ and $G=G_n(w)$ for brevity.  Let $\zeta:D\to G$ be the concatenation of the natural maps $\xi:D\to F_n$ and $\eta:F_n\to G$.  Given a  subgroup $G'$ of $G$, there is a corresponding subgroup $D'=\zeta^{-1}(G')$ of $D$.

We shall construct an Eilenberg--Mac~Lane space $X$ for $D$ as follows.  Let $\Gamma$ be a finite connected 1-complex such that $\pi_1(\Gamma)=F_n$, and realize the element $w\in F_n$ as a map $i:C\to\Gamma$ where $C$ is a circle.  Now $X$ is constructed by gluing two copies of $\Gamma$ to either end of the cylinder $C\times [0,1]$, where the gluing maps are given by $i$.

The covering space $X'$ of $X$ that corresponds to $D'$ is easy to construct.  Let $F'=\eta^{-1}(G')$ be the pre-image of $G'$ in $F_n$ and let $\Gamma'$ be the covering space of $\Gamma$ corresponding to $F'$.  Let $\{i'_j:C'_j\to\Gamma'\mid j\in J\}$ be the complete set of lifts of $i$ to $\Gamma'$.  Then $X'$ is constructed by gluing two copies of $\Gamma'$ together along cylinders $\{C'_j\times[0,1]\mid j\in J\}$, where the gluing maps of $C'_j\times[0,1]$ are copies of $i'_j$.  The resulting space $X'$ is a covering space of $X$, and it is easy to see that $\pi_1(X')=D'$.  As $X'$ is aspherical, we can use it to compute the homology of $D'$.

Now consider the natural presentation complex $Y$ for $G$, constructed by gluing a 2-cell $E$ to $\Gamma$ along $i$.  Let $Y'$ be the covering space of $Y$ with fundamental group $G'$.  It is constructed from $\Gamma'$ by attaching 2-cells $\{ E'_j\mid j\in J\}$ using the attaching maps $\{i'_j\mid j\in J\}$.  This presentation complex is not \emph{a priori} aspherical, although in the context that concerns us, when $G$ is torsion-free, it is.

\begin{thm}[Lyndon \cite{lyndon_cohomology_1950}]\label{t: Lyndon's Theorem}
If $w$ is not a proper power then $Y$ is aspherical.
\end{thm}

Even if $Y'$ is not aspherical, it is the 2-skeleton of an Eilenberg--Mac~Lane space for $G'$.  Therefore there is an epimorphism $H_2(Y')\to H_2(G')$.  Theorem \ref{t: Theorem A} is now an immediate consequence of the following lemma.

\begin{lem}\label{l: H2 of D' is H2 of Y'}
For $D'$ and $Y'$ as defined above,
\[
H_2(D')\cong H_2(Y').
\]
\end{lem}
\begin{pf}
Let $C'=\coprod_{j\in J} C'_j$ and let $i':C'\to\Gamma'$ be the map whose restriction to $C'_j$ is $i'_j$.  The Mayer--Vietoris sequence for $X'$ gives
\[
0\to H_2(X')\to H_1(C')\stackrel{j}{\to} H_1(\Gamma')\oplus H_1(\Gamma')\to\cdots
\]
where $j(x)=(i'_*(x),-i'_*(x))$.  So $H_2(X')\cong\ker j\cong\ker i'_*$.

The Mayer--Vietoris sequence applied to $Y'$ gives
\[
0\to H_2(Y')\to H_1(C')\stackrel{i'_*}{\to} H_1(\Gamma')\to\cdots
\]
so $H_2(Y')\cong\ker i'_*\cong H_2(X')\cong H_2(D')$.
\end{pf}

Concatenating the isomorphism $H_2(D')\cong H_2(Y')$ with the natural surjection $H_2(Y')\to H_2(G')$ completes the proof of Theorem \ref{t: Theorem A}, relating the homology of subgroups of the double $D$ to the homology of subgroups of the one-relator group $G$.   Note that if $w$ is not a proper power then $Y$ is aspherical by Theorem \ref{t: Lyndon's Theorem}, so $Y'$ is aspherical and the natural map $H_*(Y')\to H_*(G')$ is an isomorphism.

Corollary \ref{c: Corollary B} follows immediately from Theorem \ref{t: Theorem A} and the following lemma.

\begin{lem}\label{l: Virtual b_2}
If $G'$ is  an index-$k$ subgroup of $G$ and $D'=\zeta^{-1}(G')$ then
\[
\beta_2(D')= \beta_1(G')+k(2-n)-1
\]
\end{lem}
\begin{pf}
The Euler characteristic of the complex $Y$ is $2-n$.  If $Y'$ is the covering space of $Y$ corresponding to $G'$ then it follows from the multiplicativity of Euler characteristic that
\[
1-\beta_1(Y')+\beta_2(Y')=k(2-n).
\]
Applying Lemma \ref{l: H2 of D' is H2 of Y'} together with the fact that $\beta_1(G')=\beta_1(Y')$ and rearranging the equation, the result follows.
\end{pf}

The results of this section allow us to find subgroups of the double $D$ with positive second Betti number by pulling back subgroups of the corresponding one-relator group $G$.  The resulting subgroups of $D$ are very special---for instance, they are invariant under the natural involution of $D$ obtained by swapping the factors.  In the remainder of this section, we characterize how much is lost by restricting our attention to these subgroups.

Recall that $\xi:D\to F_n$ is the natural retraction map and $\eta:F_n\to G$ is the quotient map.  The next lemma asserts that, to answer Question \ref{q: b_2}, nothing is lost in looking at subgroups pulled back using $\xi$.

\begin{lem}\label{l: Pullbacks from the free group lose nothing}
The double $D$ has a finite-index subgroup $D'$ with $\beta_2(D')>0$ if and only if the free group $F_n$ has a finite-index subgroup $F'$ such that $\beta_2(\xi^{-1}(F'))>0$.  Furthermore, if $D'\lhd D$ and $w\in D'$ then $F'$ can also be taken to be normal in $F_n$ and to contain $\xi(w)$.
\end{lem}
\begin{pf}
One implication is immediate.  For the converse, let $D'$ be a finite-index subgroup of $D$ that corresponds to a finite-sheeted covering space $X'$ of $X$.  Let $C'$ be the preimage of the circle $C$ in $X'$ and let $\Gamma'_1$ and $\Gamma'_2$ be the preimages of the copies of the 1-complex $\Gamma$ in $X$.  Let $i'_j$ be the map $C'\to\Gamma'_j$ induced by pushing $C'$ into $\Gamma'_j$.  Cutting $X'$ along $C'$ divides it into two pieces homotopy equivalent to $\Gamma'_1$ and $\Gamma'_2$.  The Mayer--Vietoris sequence for $X'$ yields
\[
0\to H_2(X')\to H_1(C')\stackrel{j'}{\to} H_1(\Gamma'_1)\oplus H_1(\Gamma'_2)\to\cdots
\]
where $j'(x)=(i'_{1*}(x),-i'_{2*}(x))$.  If $H_2(X')\neq 0$ then it follows that $i'_{1*}$ is not injective.  Let $\bar{X}$ be obtained from two copies of $\Gamma'_1$ by gluing $C'\times[0,1]$, a disjoint union of cylinders, to each copy using $i'_1$ as the gluing map.  As before, the Mayer--Vietoris sequence gives
\[
0\to H_2(\bar{X})\to H_1(C')\stackrel{\bar{j}}{\to} H_1(\Gamma'_1)\oplus H_1(\Gamma'_1)\to\cdots
\]
where $\bar{j}(x)=(i'_{1*}(x),-i'_{1*}(x))$.  Because $i'_{1*}$ is not injective it follows that $\bar{X}$ has positive second Betti number.  Let $\hat{X}$ be a connected component of $\bar{X}$ with $\beta_2(\hat{X})> 0$ and let $\hat{\Gamma}=\Gamma'_1\cap\hat{X}$.  Set $F'=\pi_1(\hat{\Gamma})$.  Now $\pi_1(\hat{X})=\xi^{-1}(F')$, as required.

If $D'$ is a normal subgroup of $D$ then the covering map $X'\to X$ is regular.  It follows that $\Gamma'_1$ and $\hat{\Gamma}$ are both regular covering spaces of $\Gamma$.  Furthermore, if $w\in D'$ then the map $C'\to C$ restricts to a homeomorphism on each connected component.  It follows that the map $C\to\Gamma$ lifts to a map $C\to\hat{\Gamma}$, and thence that $\xi(w)\in F'=\pi_1(\hat{\Gamma})$.
\end{pf}

The assumption that a suitable finite-index subgroup $D'$ can be pulled back using $\zeta$ is more restrictive, but in a way that is easy to characterize.

\begin{lem}
The one-relator group $G$ has a finite-index subgroup $G'$ with
\[
\beta_2(\zeta^{-1}(G'))>0
\]
if and only if the double $D$ has a finite-index normal subgroup $D'$ with $w\in D'$ and $\beta_2(D')>0$.
\end{lem}
\begin{pf}
Suppose $G'$ is a finite-index subgroup of $G$ with $\beta_2(\zeta^{-1}(G'))>0$.  Let $\hat{G}$ be the normal core of $G'$ (defined to be the intersection of all the conjugates of $G'$), so $\hat{G}$ is a finite-index normal subgroup of $G$.  Therefore $D'=\zeta^{-1}(\hat{G})$ is a finite-index normal subgroup of $D$, and because $\zeta(w)=1$ we have that $w\in D'$.  It is an easy exercise with the transfer map to see that $H_2(\zeta^{-1}(G');\Q)$ embeds in $H_2(D';\Q)$ and hence that $\beta_2(D')>0$.

For the converse, let $D'$ be a finite-index normal subgroup of $D$ with $w\in D'$ and $\beta_2(D')>0$.  By Lemma \ref{l: Pullbacks from the free group lose nothing}, $F_n$ has a finite-index subgroup $F'$ such that $\beta_2(\xi^{-1}(F'))>0$.  Furthermore, Lemma \ref{l: Pullbacks from the free group lose nothing} asserts that $F'$ is a normal subgroup of $F_n$ and that $\xi(w)\in F'$.  Therefore, if $G'=\eta(F')$ then $D'=\zeta^{-1}(G')$, as required.
\end{pf}

Note that the hypothesis that $D'$ is normal is not in itself restrictive---one can always ensure this by passing to the normal core.  However, the combined hypotheses that $D'$ is normal and that $w\in D'$ are a genuine restriction on the finite-index subgroups that we consider.

\section{The Alexander polynomial}\label{s: Alexander polynomial}

In this section, we use the Alexander polynomial and Corollary \ref{c: Corollary B} to give a positive answer to Gromov's question for certain examples.

Let $G$ be a finitely presented group and $\phi:G\to\Z$ an epimorphism.  Let $Y$ be a connected complex with $\pi_1(Y)\cong G$ and let $Y_\infty\to Y$ be the $\Z$-covering with $\pi_1(Y_\infty)\cong\ker\phi$.  Then $H_1(Y_\infty)$ is a finitely presented module over $R=\Z[\Z]=\Z[t,t^{-1}]$, which depends only on $\phi$.  One can then define its elementary ideals, Alexander polynomial etc \cite{fox_free_1954}.

Specializing to the case of interest to us here, let $w\in F_n$ and let $\phi: G_n(w)\to\Z$ be an epimorphism.  By applying Nielsen transformations, we can choose a basis $x_1,\ldots, x_{n-1}, z$ for $F_n$ such that $\phi(z)=1$ and $\phi(x_i)=0$ whenever $1\leq i\leq n-1$.  In particular, the exponent sum of $z$ in $w$ is zero.  Let $Y=Y^{(1)}\cup E$ be the 2-complex associated with the one-relator presentation of $G_n(w)$, so $Y^{(1)}$ is a wedge of $n$ circles and $E$ is a 2-cell.  Let $Y_\infty^{(1)}\cup E_\infty\to Y$ be the $\Z$-covering induced by $\phi$.  Then $H_1(Y^{(1)}_\infty)$ is a free $R$-module of rank $n-1$, with basis corresponding to $x_1,\ldots,x_{n-1}$, and $H_1(\partial E_\infty)\cong R$.  Hence
\[
H_1(Y_\infty)\cong R^{n-1}/((f_1,\ldots,f_{n-1}))
\]
where $f_i\in R$ for each $1\leq i\leq {n-1}$.

\begin{defn}
The \emph{Alexander polynomial} $\Delta(t)=\Delta_\phi(t)$ is defined to be the greatest common divisor of $\{f_i\mid 1\leq i\leq n-1\}$.  It is well-defined up to multiplication by a unit $\pm t^r$ in $R$.
\end{defn}

The polynomials $f_i$ are easily determined as follows.  Since the exponent sum of $z$ in $w$ is zero, $w$ can be expressed as a word $u$ in $\{x_i^{(j)}=z^{-j}x_iz^j\mid1\leq i\leq n-1, j\in\Z\}$.  Then, for $1\leq i\leq n-1$,
\[
f_i(t)=\sum_j \sigma_i^{(j)}t^j
\]
where $\sigma_i^{(j)}$ is the exponent sum of $x_i^{(j)}$ in $u$.

The above discussion may be carried out with the coefficients $\Z$ replaced by $\F_p$ or $\C$, for example.

\begin{thm}\label{t: Roots of unity}
Let $w\in F_n$.  If there is an epimorphism $\phi:G_n(w)\to\Z$ such that $\Delta_\phi(t)$ has a root $\omega$ with $\omega^k=1$ for some $k$ then there is a finite-index subgroup $G'$ with $\beta_1(G')>1+k(n-2)$.  Therefore $D_n(w)$ contains a surface subgroup.
\end{thm}

To prove the theorem we consider the finite-index subgroups $G_k=\phi^{-1}(k\Z)$ of $G_n(w)$.  The theorem follows immediately from the combination of Corollary \ref{c: Corollary B} and Lemma \ref{l: Alexander polynomial}.  In the context of the homology of the branched cyclic covers of links, the lemma is essentially contained in \cite{sumners_homology_1974}.  (For knots the corresponding result is due to Goeritz \cite{goeritz_die_1934}.)  We follow closely the elegant treatment of Sumners.

\begin{lem}\label{l: Alexander polynomial}
For any $k$, $\beta_1(G_k)>k(n-2)+1$ if and only if some root of $\Delta_\phi(t)$ is a $k$th root of unity.
\end{lem}
\begin{pf}
Let $Y_k$ be the finite-sheeted covering space of $Y$ associated to $G_k$.  We consider homology with coefficients in $\C$.  Let $S=R\otimes\C=\C[\Z]=\C[t,t^{-1}]$.  Because $k\Z=\langle t^k\rangle$ acts on $Y_\infty$ with quotient $Y_k$, there is a short exact sequence of chain complexes
\[
0\to C_*(Y_\infty;\C)\stackrel{t^k-1}{\to}C_*(Y_\infty;\C)\to C_*(Y_k;\C)\to 0
\]
that induces a long exact sequence in homology.
\[
\cdots\to H_1(Y_\infty;\C)\stackrel{t^k-1}{\to}H_1(Y_\infty;\C)\to H_1(Y_k;\C)\to H_0(Y_\infty;\C)\stackrel{t^k-1}{\to}H_0(Y_\infty;\C)
\]
But $\dim_\C H_0(Y_\infty;\C)=1$ and $t^k-1$ acts as $0$ on $H_0(Y_\infty;\C)$ so
\begin{equation}\label{eq: Roots of unity}
\dim_\C H_1(Y_k;\C)=1+\dim_\C \mathrm{coker}(t^k-1).
\end{equation}
Because $S$ is a principal ideal domain, it is easy to see that
\[
H_1(Y_\infty;\C)\cong S^{n-2}\oplus S/(\Delta(t))
\]
where $\Delta(t)=\Delta(t)\otimes 1\subset S$.  Clearly
\[
\dim_\C(\mathrm{coker}(t^k-1:S\to S))=k.
\]
To calculate $\dim_\C(\mathrm{coker}(t^k-1:S/(\Delta(t))\to S/(\Delta(t))))$, consider a cyclic summand
\[
V=S/(t-\alpha)^m
\]
of $S/(\Delta(t))$, and the restriction of $t^k-1$ to $V$.  We have
\[
\mathrm{coker}(t^k-1:V\to V)=S/((t-\alpha)^m, t^k-1).
\]
Over $\C$, $t^k-1$ factors as
\[
t^k-1=\prod_j(t-\omega_j)
\]
where $\omega_1,\ldots,\omega_k$ are the $k$th roots of unity. If $\alpha=\omega_j$ for some $j$ then the highest common factor of $(t-\alpha)^m$ and $t^k-1$ is $t-\alpha$, in which case the cokernel is $S/(t-\alpha)\cong\C$.  Otherwise $(t-\alpha)^m$ and $t-\omega_j$ are relatively prime, so the cokernel is trivial.

Therefore, the dimension of the cokernel of $t^k-1$ on $S/(\Delta(t))$ is equal to the number, $r$ say, of $\alpha$'s that are $k$th roots of unity.  By (\ref{eq: Roots of unity}) above we then have that
\[
\beta_1(G_k)=1+k(n-2)+r
\]
and the result follows.
\end{pf}

We can also exploit work of Howie \cite{howie_free_1998} to provide more examples using only the Alexander polynomial modulo a prime.  The case $n=2$ is particularly simple.  First, in this case $H_1(Y_\infty)$ is the cyclic module $R/(\Delta(t))$ and, secondly, if $w\notin [F_2,F_2]$ (which we may assume, otherwise $\beta_1(G_2(w))=2$) then $\phi$ is unique (up to automorphisms of $\Z$).  Hence we shall assume that $w\in F_2\smallsetminus [F_2,F_2]$ and we may write $\Delta_w(t)=\Delta_\phi(t)$.

\begin{thm}\label{t: Alexander polynomial mod p}
Let $w\in F_2\smallsetminus [F_2,F_2]$.  If $\Delta_w(t)\equiv 0$ mod $p$ for some prime $p$ then $G_2(w)$ is large. Therefore, $D_2(w)$ contains a surface subgroup.
\end{thm}
\begin{pf}
Recall that $H_1(Y_\infty)\cong R/(\Delta(t))$ (writing $\Delta(t)=\Delta_w(t)$).  Hence
\[
H_1(Y_\infty;\F_p)\cong H_1(Y_\infty)\otimes \F_p\cong (R\otimes\F_p)/(\Delta(t)\otimes 1).
\]
If $\Delta(t)\equiv 0$ mod $p$ then $H_1(Y_\infty;\F_p)\cong R\otimes\F_p\cong\F_p[t,t^{-1}]$.  As pointed out in \cite{button_largeness_2008}, the argument of Proposition 2.1 and Theorem 2.2 of \cite{howie_free_1998} then shows that $G_2(w)$ is large.
\end{pf}

\section{Positive $C'(1/6)$ one-relator groups}\label{s: C'(1/6)}

We refer the reader to \cite{lyndon_combinatorial_1977} for the definition of the small-cancellation condition $C'(\lambda)$.  In this section we shall exploit a theorem of Wise \cite{wise_residual_2001} to deduce that doubles of $F_2$ along positive $C'(1/6)$ words have surface subgroups.

\begin{defn}
A graph of free groups is \emph{clean} if it can be realized as a graph of spaces in which every vertex space is a connected 1-complex and every edge map is injective.  Furthermore, we shall call such a graph of groups \emph{properly clean} if no edge map is $\pi_1$-surjective.
\end{defn}

The following theorem is essentially due to Wise, who proved that certain one-relator groups have a finite-index subgroup that is the fundamental group of a clean graph of groups.  We shall briefly explain why this graph of groups can be taken to be not just clean, but properly clean.  Recall that a subgroup $H$ of a group $G$ is \emph{malnormal} if, for all $g\in G$, $gHg^{-1}\cap H\neq 1$ implies that $g\in H$.

\begin{thm}[Wise \cite{wise_residual_2001}]\label{t: Wise 1-relator}
If $w\in F_n$ is a positive $C'(1/6)$ word then the one-relator group $G_n(w)$ has a finite-index subgroup that is the fundamental group of a finite, properly clean graph of free groups with finitely generated edge groups.
\end{thm}
\begin{pf}
Let $G=G_n(w)$, where $w$ is a positive $C'(1/6)$ word.  Theorem 1.1 of \cite{wise_residual_2001} asserts that $G$ has a finite-index subgroup that splits as an amalgamated product $G'=A*_C B$, where $A$ and $B$ are free and $C$ is malnormal in each of $A$ and $B$.  If $C$ is equal to either $A$ or $B$ then the amalgamated product decomposition is trivial and $G'$ is free.   Therefore, we may assume that $C$ is a proper, malnormal subgroup of both $A$ and $B$.

Now, it follows from Theorem 11.3 of \cite{wise_residual_2002} that $G'$ has a finite-index subgroup $\hat{G}$ for which the induced graph-of-groups decomposition is clean.  Furthermore, as $C$ is a proper, malnormal (and hence infinite-index) subgroup of $A$ and $B$, the edge maps of the decomposition of $\hat{G}$ are never $\pi_1$-surjective.  Therefore, the graph of groups for $\hat{G}$ is properly clean.
\end{pf}

Theorem \ref{t: Theorem C} now follows from the following result.  Note that the hypothesis that the graph of spaces is \emph{properly} clean makes things much easier; it is a long-standing question whether all free-by-cyclic group have virtual first Betti number greater than $2$.

\begin{thm}
If $G$ is the fundamental group of a finite, properly clean graph of free groups with finitely generated edge groups then either $G$ is cyclic or $G$ surjects a non-abelian free group.
\end{thm}
\begin{pf}
By hypothesis, $G$ is the fundamental group of a finite graph of spaces $X$, in which every vertex space is a connected 1-complex, every edge space is a finite connected 1-complex and every edge map is injective but not $\pi_1$-surjective.  Let $e$ be an edge incident at distinct vertices $u$ and $v$.  Denote the corresponding edge and vertex spaces by $X_e$, $X_u$ and $X_v$ in the obvious manner.  The subspace of $X$ consisting of the union of $X_u$, $X_v$ and the cylinder $X_e\times I$ is homotopy equivalent to a 1-complex, in such a way that both $X_u$ and $X_v$ embed as subcomplexes.  Without loss of generality, therefore, the underlying graph of $X$ can be taken to be a wedge of $r$ circles.

There is a natural surjection from $G$ to the fundamental group of the underlying graph.  Therefore, if $r>1$ then the result follows.  Furthermore, if $r=0$ then $G$ is free.  The case in which $r=1$ remains.  We are therefore reduced to the case in which $X$ has a single vertex space $\Gamma=X_v$, with two finite isomorphic subcomplexes $\Gamma_1$ and $\Gamma_2$ which are the images of the two edge maps.

Suppose $\langle \pi_1(\Gamma_1),\pi_1(\Gamma_2)\rangle$ is strictly smaller than $\pi_1(\Gamma)$.  Then there is a non-separating 1-cell $\epsilon$ of $\Gamma$ in the complement of $\Gamma_1\cup\Gamma_2$.  Let $X'$ be the subspace of $X$ obtained by deleting $\epsilon$.  By van Kampen's Theorem,
\[
G=\pi_1(X')*\Z.
\]
The subspace $X'$ still has the structure of a graph of spaces, the underlying graph of which is a circle, so $\pi_1(X')$ surjects $\Z$.  Therefore $G$ surjects a non-abelian free group.

Finally, we are left with the case in which $\langle \pi_1(\Gamma_1),\pi_1(\Gamma_2)\rangle=\pi_1(\Gamma)$.  In this case, there is a 1-cell $\epsilon\in\Gamma_1\smallsetminus\Gamma_2$ (otherwise, $\pi_1(\Gamma_2)=\pi_1(\Gamma)$).  The edge space of $X$ contains a rectangle of the form $\epsilon\times I$.  Because this rectangle has a free edge, it can be collapsed, yielding a new graph of spaces $\hat{X}$ whose edge space has one fewer 1-cells.  The result now follows by induction.
\end{pf}

\section{Baumslag--Solitar groups}\label{s: Baumslag-Solitar groups}

One might naively hope to use these methods to find a surface subgroup in every one-ended double $D_2(w)$---that is, one might hope that every freely indecomposable two-generator, one-relator group might have virtual first Betti number greater than 1.  This is too optimistic.

\emph{Baumslag--Solitar groups} are two-generator, one-relator groups with presentations
\[
BS(p,q)\cong\langle a,b\mid b^{-1}a^pb=a^q\rangle.
\]
Note that the abelianization of $BS(p,q)$ is $\Z\oplus(\Z/(p-q))$.  In particular, there is always a homomorphism onto $\Z$ that sends $b$ to 1 and $a$ to $0$.  However, there is a sharp dichotomy in the behaviour of the virtual Betti numbers of these groups.

Suppose that $p$ and $q$ have a common factor $k>1$, so $p=kp'$ and $q=kq'$.  Then $BS(p,q)$ is obtained by adjoining a $k$th root to $BS(p',q')$, and admits a map $BS(p,q)\to\Z/k\Z$ that sends $b$ to $0$ and $a$ to $1$.  The kernel of this map is an index-$k$ subgroup $K$ with presentation
\[
K\cong\langle \alpha,\beta_1,\ldots,\beta_k\mid\beta_1^{-1}\alpha^{p'}\beta_1=\alpha^{q'},\ldots,\beta_k^{-1}\alpha^{p'}\beta_k=\alpha^{q'}\rangle
\]
where $\alpha=a^k$ and $\beta_i=a^{1-i}ba^{i-1}$ in the original group.  There is a homomorphism from $K$ onto a free group of rank $k$ given by setting $\alpha$ to be trivial. Therefore in this case $BS(p,q)$ is large.

By contrast, if $p$ and $q$ are relatively prime and not both $\pm1$ then it is a theorem of Edjvet and Pride that the virtual first Betti number of $BS(p,q)$ is 1.  We shall give a proof of this theorem.  First, we recall a little of the theory of graphs of groups.

Let $G$ be the fundamental group of a graph of groups $\Gamma$.  There is an associated Bass--Serre tree $T$ with an action of $G$.  The set of vertices is in bijection with the coset space $G_v\backslash G$ and the set of edges is in bijection with the coset space $G_e\backslash G$, with adjacency given by inclusion.  The action of $G$ is by right multiplication.  See \cite{serre_arbres_1977} for details.

For a subgroup $K$ of $G$ there is an induced graph-of-groups decomposition $\Gamma^K=T/K$.  The set of vertices of $\Gamma^K$ is in bijection with the set of double cosets $G_v\backslash G/K$, and similarly the set of edges is in bijection with the set of double cosets $G_e\backslash G/K$.

For a vertex $v'$ of $\Gamma^K$ corresponding to a double coset $G_vgK$, the corresponding vertex group of $\Gamma^K$ is $K\cap G_v^g$.  (Note that this is well-defined up to conjugation in $K$.) The \emph{index} of $v'$ is defined to be the index $|G_v:G_{v'}|$.  Similarly, for an edge $e'$ of $\Gamma^K$ corresponding to a double coset $G_egK$, the corresponding edge group of $\Gamma^K$ is $K\cap G_e^g$.  The \emph{index} of $e'$ is defined to be the index $|G_e:G_{e'}|$.  The set of edges of $\Gamma^K$ incident at $v'$ that lie above $e$ is in bijection with the set of double cosets $G_e\backslash G_v/G_{v'}$.

\begin{thm}[Edjvet and Pride \cite{edjvet_concept_1984}]\label{t: Edjvet-Pride}
If $p,q$ are relatively prime and not both equal to $\pm 1$ and $K$ is a finite-index subgroup of $BS(p,q)$ then $\beta_1(K)=1$.
\end{thm}
\begin{pf}
Think of $G=BS(p,q)$ as the fundamental group of a graph of groups $\Gamma$ with one vertex $v$ with cyclic vertex group $G_v=\langle a\rangle$ and one edge $e$ with cyclic edge group $G_e$.  Without loss of generality, we can take $K$ to be normal.  We shall start by proving that the underlying graph of $\Gamma^K$ is topologically a circle.

There is an orientation on $\Gamma$ determined by the requirement that the corresponding stable letter $b$ conjugates $a^p$ to $a^q$, and this induces an orientation on $\Gamma^K$.  So it makes sense to say that an edge of $\Gamma^K$ points towards or away from a vertex.  The quotient group $K\backslash G$ acts on $\Gamma^K$ with quotient $\Gamma$, so there is one orbit of edges and one orbit of vertices.  In particular, all vertices of $\Gamma^K$ have the same index, and similarly all the edges have the same index.  Fix a vertex $v'$ of $\Gamma^K$, of index $k$.  Without loss of generality, $G_{v'}=\langle a^k\rangle$.  The set of edges incident at $v'$ that point towards $v'$ is in bijection with the set of double cosets
\[
G_e\backslash G_v/G_{v'}=\langle a^q\rangle\backslash\langle a\rangle/\langle a^k\rangle\cong\Z/(q,k)
\]
and each such edge has index $k/(n,k)$.  Similarly, the set of edges incident at $v'$ that point away from $v'$ is in bijection with the set of double cosets
\[
G_e\backslash G_v/G_{v'}=\langle a^p\rangle\backslash\langle a\rangle/\langle a^k\rangle\cong\Z/(p,k)
\]
and each such edge has index $k/(p,k)$.  But all edges have the same index, so $(p,k)=(q,k)$.  Because $p$ and $q$ are relatively prime, $(p,k)=(q,k)=1$ and there is just one incident edge pointing towards $v'$ and likewise just one incident edge pointing away from $v'$.

Therefore $\Gamma^K$ is topologically a circle. Furthermore, every vertex and every edge of $\Gamma^K$ has degree $k$, where $k$ is relatively prime to both $p$ and $q$.  Let $l$ be the number of vertices of $\Gamma^K$.  Then $K$ admits the following presentation:
\[
K\cong\langle \alpha_1,\ldots,\alpha_l,\beta\mid \alpha_1^p=\alpha_2^q,\alpha_2^p=\alpha_3^q,\ldots,\alpha_{l-1}^p=\alpha_l^q,\beta^{-1}\alpha_l^p\beta=\alpha_1^q \rangle.
\]
For each $i$ it is easy to check that
\[
\beta^{-1}\alpha_i^{p^l}\beta=\alpha_i^{q^l}
\]
and so $\alpha_i$ has finite order in the abelianization of $K$ unless $p$ and $q$ are both $\pm 1$.  Therefore $\beta_1(K)=1$ as required.
\end{pf}

\section{Virtually geometric words}\label{s: Virtually geometric}

As we saw in Section \ref{s: Baumslag-Solitar groups}, the methods of Section \ref{s: Mayer-Vietoris} cannot answer Gromov's question for $D_2(w)$ when
\[
w=b^{-1}a^pba^q\in\langle a,b\rangle
\]
say.  In this section, we develop an alternative approach to finding surface subgroups of doubles.  The idea is simple: if a one-ended double has a finite-index subgroup that is the fundamental group of a compact 3-manifold then it will contain a surface subgroup, coming from the boundary.

\begin{defn}
A finite subset $S$ of a free group $F$ is called \emph{geometric} if $F$ can be realized as the fundamental group of a handlebody $V$ in such a way that a set of loops realizing $S$ is freely homotopic to an embedded multicurve in the boundary of $V$.  We say that $w\in F$ is geometric if $\{w\}$ is geometric.
\end{defn}

If $w$ is geometric then it follows that $G_n(w)$ and $D_n(w)$ are both the fundamental groups of compact 3-manifolds.  (For $G_n(w)$, simply glue on a thickened 2-cell.  For $D_n(w)$, double $V$ along thickened cylinders.)  Therefore, there are examples of non-geometric $w$.  For instance, the Baumslag--Solitar group $BS(p,q)$ is never a 3-manifold group unless $p=\pm q$ or $pq=0$ \cite{heil_finitely_1975,jaco_roots_1975} so $w=b^{-1}a^pba^q$ is not geometric.

\begin{defn}
An element $w\in F$ is \emph{virtually geometric} if $F$ has a finite-index subgroup $F'$ such that the set
\[
\{g^{-1}w^{n_g} g\mid gF'\in F/F'\}
\]
is a geometric subset of $F'$, where $g$ ranges over a set of coset representatives of $F/F'$ and $n_g$ is the minimal positive integer such that $g^{-1}w^{n_g}g\in F'$ for each $g$.
\end{defn}

If $w$ is virtually geometric then $D_n(w)$ has a finite-index subgroup $D'$ that is the fundamental group of a compact 3-manifold with boundary.  We therefore have the following.

\begin{lem}\label{l: Geometric implies surface}
If $w\in F_n$ is virtually geometric and $D_n(w)$ is one-ended then $D_n(w)$ has a surface subgroup.
\end{lem}
\begin{pf}
Let $V'$ be a handlebody representing $F'$ and let $\gamma$ be an embedded multicurve in $\partial V'=\Sigma$ representing $\{g^{-1}w^{n_g} g\mid gF'\in F/F'\}$.  Let $N(\gamma)$ be a small neighbourhood of $\gamma$ in $\Sigma$.  Attach the two ends of $N(\gamma)\times [0,1]$ to two copies of $V'$ in the natural way.  Call the resulting 3-manifold $M$.  Then $\pi_1(M)\cong D'=\xi^{-1}(F')$.

The boundary of $M$ is homeomorphic to the double of $\Sigma\smallsetminus\mathrm{Int}(N(\gamma))$.  As no component of $\gamma$ bounds a disc in $\Sigma$, every component $B$ of $\partial M$ has non-positive Euler characteristic.  If the map $\pi_1(B)\to\pi_1(M)$ is not injective then it follows from the Loop Theorem and Dehn's Lemma that $D'$, and hence $D_n(w)$, is not one-ended.
\end{pf}

This gives another approach to finding surface subgroups of doubles, which we shall apply to the Baumslag--Solitar case.

\begin{thm}\label{t: BS are virtually geoemtric}
If $(p,q)=1$ then the Baumslag--Solitar word $w=b^{-1}a^pba^q$ is virtually geometric.  Therefore, if $D_2(w)$ is one-ended then it contains a surface subgroup.
\end{thm}

\begin{figure}[htp]
\begin{center}
\includegraphics[width=0.7\textwidth]{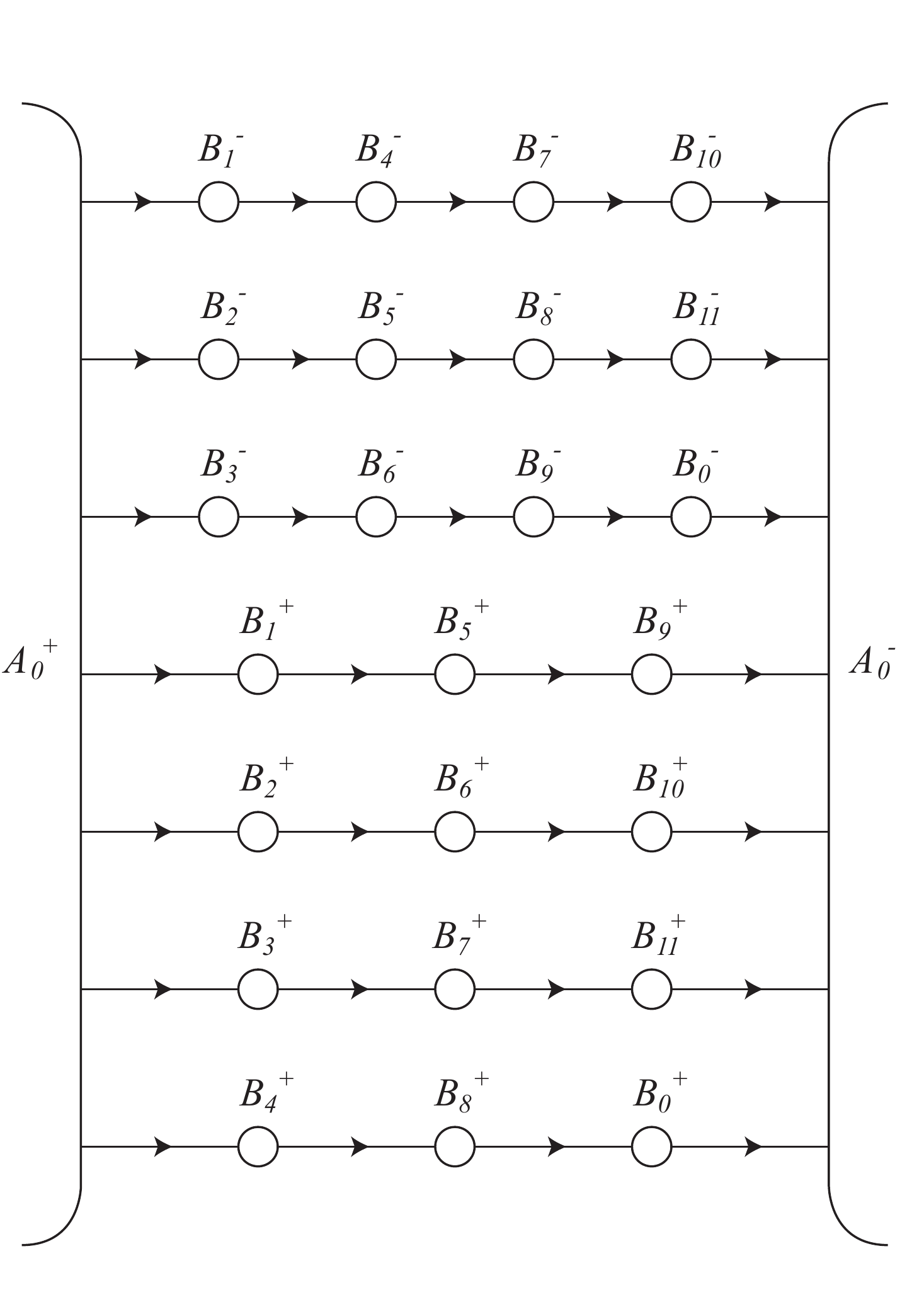}
\caption{The Heegaard diagram in the case $p=3$, $q=4$.}\label{Pre-cover figure}
\end{center}
\end{figure}

\begin{pf}
First assume that $p,q>0$.  Let $F=F_2=\langle a,b\rangle$ and let $\psi:F\to\Z/pq$ be the epimorphism defined by $\psi(a)=1$ and $\psi(b)=0$.  Let $w=b^{-1}a^pba^q\in F$.  Since $(p,q)=1$ we have $(p+q,pq)=1$ and hence the smallest positive $n$ such that $w^n\in\ker\psi$ is $n=pq$.  Let $\tilde{F}=\ker\psi$ and let $\tilde{w}=w^n\in\tilde{F}$.

Let $X$ be the wedge of two oriented circles $\alpha$ and $\beta$ representing $a$ and $b$ respectively, so $\pi_1(X)\cong F$.  Let $\rho:\tilde{X}\to X$ be the $\Z/n$-covering defined by $\psi$.  More precisely, let $x\in X$ be the wedge-point and let $\rho^{-1}(x)=\{x_i\mid i\in\Z/n\}\subset\tilde{X}$.  Then $\rho^{-1}(\alpha)$ consists of $n$ arcs $\alpha_i$ where $\alpha_i$ joins $x_i$ to $x_{i+1}$, and $\rho^{-1}(\beta)$ consists of $n$ loops $\beta_i$, where $\beta_i$ is based at $x_i$, for each $i\in\Z/n$.  Taking $\alpha_1\cup\ldots\cup\alpha_{n-1}$ as a maximal tree for $\tilde{X}$, we see that $\tilde{F}$ has basis $b_0,\ldots, b_{n-1}, a_0$, where $b_i=[\beta_i]$ and $a_0=[\alpha_0]$.

Consider $\tilde{w}\in\tilde{F}$.  Lifting the loop in $X$ representing $w^n$ to $\tilde{X}$ we see that successive occurrences of $b_i^{\pm 1}$'s in $\tilde{w}$ are of the form $b_i^{-1}$, $b_{i+p}$ and $b_j$, $b_{j+q}^{-1}$, and that $a_0$ occurs in a syllable $b_i^{-1}a_0b_{i+p}$ (respectively $b_ja_0b_{j+q}^{-1}$) if and only if $0\in [i,i+p)$ (respectively $0\in[j,j+q)$).  (Here, the intervals are to be interpreted as cyclic intervals in $\Z/n$ in the obvious way.)

Let $\tilde{V}$ be a handlebody of genus $n+1$, with a complete set of meridian discs $B_0,\ldots,B_{n-1},A_0$ corresponding to the basis $b_0,\ldots,b_{n-1},a_0$ for $\tilde{F}=\pi_1(\tilde{V})$.  Cutting $\tilde{V}$ along these discs gives a 3-ball $E$ with $n+1$ pairs of discs $B_0^+,B_0^-,\ldots,B_{n-1}^+,B_{n-1}^-,A_0^+,A_0^-$ in $\partial E$.  A simple loop on $\partial\tilde{V}$ gives rise to a Heegaard diagram of disjoint arcs in $\partial E$ with their endpoints on the boundaries of these discs; we shall show that $\tilde{w}$ is represented by a simple loop on $\partial\tilde{V}$ by constructing the corresponding Heegaard diagram.

Let $A_0^+,A_0^-$ be disjoint discs on $\partial E$.  Consider the $p$ orbits of the map $i\mapsto i+p$ in $\Z/n$; these may be labelled $\Omega_k$ for $k=1,2,\ldots,p$, where $k\in\Omega_k$.  Similarly, let $\Lambda_l$, for $l=1,2,\ldots,q$, be the orbits of the map $j\mapsto j+q$ in $\Z/n$, where $l\in\Lambda_l$.  Note that $0\in[k-p,k)$ and $0\in[l-q,l)$.

For each $\Omega_k$, draw $q$ disjoint discs on $\partial E$ labelled $B_i^-$, $i\in\Omega_k$, together with $q-1$ disjoint oriented arcs from $B_i^-$ to $B_{i+p}^-$, $i\in\Omega_k\smallsetminus\{k-p\}$.  Attach each of these linear graphs to $A_0^+$ and $A_0^-$ by inserting an oriented arc from $A_0^+$ to $B_k^-$ and from $B_{k-p}^-$ to $A_0^-$.  Similarly, for each orbit $\Lambda_l$, draw $p$ disjoint disks labelled $B_j^+$, $j\in\Lambda_l$, and $p-1$ arcs from $B_j^+$ to $B_{j+q}^+$, $j\in\Lambda_l\smallsetminus\{l-q\}$.  Again, attach each of these resulting linear graphs to $A_0^+$ and $A_0^-$ by adding an arc from $A_0^+$ to $B_l^+$ and from $B_{l-q}^+$ to $A_0^-$.  Figure \ref{Pre-cover figure} illustrates the case $p=3, q=4$.

Finally, to obtain $\tilde{V}$, identify $B_i^+$ with $B_i^-$ so that the endpoint of the incoming (respectively outgoing) arc at $B_i^+$ is identified with the endpoint of the outgoing (respectively incoming) arc at $B_i^-$, for each $i\in\Z/n$, and identify $A_0^+$ and $A_0^-$ so that the two endpoints of each of the $p+q$ linear graphs described above are identified.  Then the arcs described define an oriented simple loop on $\partial\tilde{V}$ that represents $\tilde{w}\in\pi_1(\tilde{V})\cong\tilde{F}$.

The case in which $p$ and $q$ have opposite sign is completely analogous; we leave the details to the reader.
\end{pf}

It seems difficult to prove that a given word $w$ is not virtually geometric.

\begin{qu}\label{q: Virtually geometric}
Is every element $w\in F_n$ virtually geometric?
\end{qu}

By Lemma \ref{l: Geometric implies surface}, a positive answer to Question \ref{q: Virtually geometric} would imply a positive answer to Gromov's question for all doubles $D_n(w)$.  One word is particularly notable for eluding all the techniques of this paper.  Baumslag \cite{baumslag_non-cyclic_1969} proved that every finite quotient of the non-cyclic one-relator group
\[
G=\langle a,b\mid a^2=b^{-1}a^{-1}bab^{-1}ab\rangle
\]
is cyclic.  The group $G$ therefore has virtual first Betti number equal to one, and the relation $w=a^{-2}b^{-1}a^{-1}bab^{-1}ab$ is not a geometric element of $F_2=\langle a,b\rangle$.  (If it were, then $G$ would be a 3-manifold group and hence residually finite.)

\begin{qu}
Is $w=a^{-2}b^{-1}a^{-1}bab^{-1}ab$ virtually geometric?  Does $D_2(w)$ contain a surface subgroup?
\end{qu}

\bibliographystyle{plain}

\end{document}